\documentclass[11pt]{article}

\usepackage{graphicx}

\usepackage[french]{babel}
\usepackage[utf8]{inputenc}

\usepackage{amssymb}
\def\E{{\mathbb{E}}}
\def\P{{\mathbb{P}}}

\def\R{{\mathbb{R}}}

\newtheorem{theo}{Th\'eor\`eme}
\newtheorem{cor}[theo]{Corollaire}
\newtheorem{prop}[theo]{Proposition}

\newtheorem{defin}[theo]{D\'efinition}

\newenvironment{dem}[1][Proof]{{\it Démonstration : } \ }{\rule{0.5em}{0.5em}}

\def\tend#1{\mathop{\hbox to #1{\rightarrowfill}}\limits}

\title {Processus de Lévy avec changements de rythmes}

\author{Christiane {\sc Cocozza-Thivent}\\
$ $ \\
anciennement membre du Laboratoire d'Analyse et de Mathématiques Appliquées, \\UMR CNRS 8050\\
Université Paris-Est Marne-la-Vallée \\
$ $ \\
adresse électronique : cocozza.christiane@orange.fr
}

\date{}

\begin{document}

\maketitle

\vspace{1cm}

\noindent

\begin{center}
Abstract
\end{center}

This paper introduces Switching Processes, called SP. Their constructions are inspired by the PDMP's ones (PDMP stands for Piecewise Deterministic Markov Process). A Markov process, called the intrinsic process, replaces the PDMP's flow. Jumps are added ; they occur randomly as their locations ; their distributions depend on the process's trajectory between them. When the intrinsic process is a Levy process, thanks to its Lévy-Itô decomposition as a semi-martingale, we obtain the expected Kolmogorov equations for the SP. The results are extended to Itô-Lévy processes, in particular to diffusion processes.

\begin{center}
Résumé
\end{center}

Dans cet article, nous introduisons les processus avec changements de rythmes. Leur construction est inspirée par celle des PDMP (Piecewise Deterministic Markov Process). Ces processus, notés SP pour Switching Processes, sont construits à partir d'un processus dit intrinsèque qui remplace le flot déterministe de la construction des PDMP. Des sauts sont ajoutés à ce processus intrinsèque. Ils se produisent à des instants aléatoires, les lois de ces instants et leurs localisations dépendent de la trajectoire du processus entre ceux-ci. Lorsque le processus intrinsèque est un processus de Lévy, son écriture comme semi-martingale (décomposition de Lévy-Itô) nous permet d'obtenir les équations de Kolmogorov auxquelles on s'attend pour le SP. Les résultats s'étendent aux processus d'Itô-Lévy et en particulier aux diffusions. 

\section{Construction des processus avec changements de ryth\-mes}

Soit $\zeta=(\zeta(t))_{t \ge 0}$ un processus càd-làg, appelé processus intrinsèque, que nous supposons à valeurs dans $\R^d$ pour simplifier. Notons ${\cal P}(\R^d)$ l'ensemble des probabilités sur $\R^d$ muni de la tribu borélienne. Les changements de rythmes consistent à lui ajouter des sauts selon un taux et un lieu de saut qui dépendent de l'état du processus.
Le taux est caractérisé par une fonction $\lambda : \R^d \to \R_+$ et le lieu du saut par une probabilité de transition $Q : \R^d \to {\cal P}(\R^d)$.
On suppose que presque sûrement quelle que soit la loi initiale de $\zeta$ il existe $\varepsilon>0$ tel que $\int_0^{\varepsilon}\lambda(\zeta(v))\, dv <+\infty$ (propriété que doit vérifier un taux de hasard).
Pour simplifier on suppose également que presque sûrement quelle que soit la loi initiale de $\zeta$ on a  $\int_0^{ +\infty} \lambda(\zeta(v))\, dv =+\infty$, ce qui entraine que pour $M$ défini ci-dessous par (\ref{defM}) on a $M(\zeta; \R^d \times \R_+)=1$.

\medskip
La terminologie "changements de rythmes" (switching en anglais) provient de certaines applications. Le modèle que nous allons présenter permet par exemple de modéliser des phénomènes dans lesquels les paramètres changent à certains instants aléatoires. Dans ce cas le processus $\zeta$ est en fait une famille de processus $(\zeta_i)_{i \in I}$ et après un saut l'évolution du phénomène qui était décrite par le processus $\zeta_i$ devient régie par le processus $\zeta_j$, $j$ étant choisi selon une probabilité $Q$ qui dépend de $i$ et de l'état du processus à l'instant du saut ; des exemples sont donnés dans le chapitre 8 de \cite{CC}. Mais, comme indiqué ci-dessus, pour simplifier nous supposons ici que $\zeta$ est à valeurs dans $\R^d$.

Dans \cite{Bect}, J. Bect s'intéresse aux processus de Markov diffusifs par morceaux, cela correspond au cas où le processus intrinsèque est une diffusion. 
De nombreux exemples d'applications sont donnés dans son introduction.

\medskip
 Posons
 \begin{equation}\label{defM}
  M(\zeta; dx,dv) = \lambda(\zeta(v))\, e^{-\int_0^v \lambda(\zeta(w))\, dw}\, Q(\zeta(v); dx)\, dv.
  \end{equation}
  
 \begin{defin}\label{defSP}
Un processus  à changements de rythmes (Switching Process ou SP)   associé à $\zeta$ et $M$ d'état initial $x_0\in \R^d$ est un processus $X=(X_t)_{t \ge 0}$ qui peut être construit de la manière suivante :
\begin{enumerate}
\item soit $\zeta^{(1)}$ un processus dont la loi est la loi de $\zeta$ sachant $\zeta(0)=x_0$,
\item la loi de $(Y_1,T_1)$ sachant $\zeta^{(1)}$ est $M(\zeta^{(1)};dx,dv)$,
\item pour $t<T_1,$ $X_t=\zeta^{(1)}(t)$ et $X_{T_1}=Y_1,$\\
$ $ \\
on suppose construits $\zeta^{(1)}, \ldots, \zeta^{(n)}, Y_1, T_1, \ldots, Y_n, T_n$ $(n \ge 1)$,
\item soit $\zeta^{(n+1)}$ un processus dont la loi sachant $\zeta^{(1)}, \ldots, \zeta^{(n)}, Y_1, T_1, \ldots, Y_n, T_n$ est la loi de $\zeta$ sachant $\zeta(0)=Y_n$,
\item la loi de $(Y_{n+1}, T_{n+1}-T_n)$ sachant $\zeta^{(1)}, \ldots, \zeta^{(n+1)}, Y_1, T_1, \ldots, Y_n, T_n$ est \\$M(\zeta^{(n+1)};dz,dv)$,
\item si $T_n \le t < T_{n+1}$, $X_t=\zeta^{(n+1)}(t-T_n),$ et $X_{T_{n+1}} = Y_{n+1}$.
\end{enumerate}

On pose $T_0=0, Y_0=x_0.$
\end{defin}

On suppose que $\lim_{n \to +\infty} T_n=+\infty,$ ce qui est le cas sous l'hypothèse $\lambda$ bornée que nous ferons ultérieurement.

\section{Approche semi-régénérative}

On définit le noyau de renouvellement $N$ sur $\R^d$ par 
\begin{equation}\label{renouvMulti}
 N(x,dy,dv)= \E(M(\zeta;dy,dv)\, / \, \zeta(0)=x),
 \end{equation}
au sens où
$$ \int_{\R^d \times \R_+} \varphi(y,v)\, N(x,dy,dv) = \E\biggl( \int_{\R^d \times \R_+} \varphi(y,v)\, M(\zeta;dy,dv) \, / \, \zeta(0)=x \biggr)$$
pour toute fonction $\varphi$ mesurable positive définie sur $\R^d \times \R_+$.

Le processus $(Y_n,T_n)_{n \ge 1}$ est un processus de renouvellement markovien de noyau de renouvellement $N$ et la loi de $(Y_1,T_1)$ sachant $X_0=x_0$ est $N(x_0, \cdot, \cdot)$. Le processus $(X_t)_{ t \ge 0}$ est un processus semi-régénératif associé au processus de renouvellement markovien $(Y_n, T_n)_{n \ge 1}$.

Nous définissons $(Z_t, A_t)$ par   $Z_t=Y_n$ et $A_t=t-T_n$ sur $T_n \le t < T_{n+1}$ $(n \ge 0)$. Nous appelons $(Z_t, A_t)_{t \ge 0}$ le CSMP (Completed Semi-Markov Process) sous-jacent au SP $(X_t)_{t \ge 0}$.

\begin{prop}\label{loiMulti}
Soit $(X_t)_{t \ge 0}$ un SP associé à $\zeta$ et $M$ donné par (\ref{defM}), $(Z_t, A_t)_{t \ge 0}$ son CSMP sous-jacent et $(Y_n, T_n)_{n \ge 1}$ le processus de renouvellement markovien associé. On pose $N_t=\sum_{n \ge 1} 1_{\{T_n \le t\}}$.

Soit $g$ une fonction mesurable positive définie sur $\R_+ \times \R^d$ et $t>0$. Alors :
$$\E(g(t,X_t)\, / \, T_{N_t}, Z_t, A_t) = \psi(T_{N_t},Z_t, A_t)\, \quad p.s.$$
avec 
$$ \psi(s,z,v) = \frac{\E\biggl( g(s+v,\zeta(v))\,  e^{-\int_0^v \lambda(\zeta(w))\, dw}  \, / \, \zeta(0)=z\biggr)}{\E\biggl( e^{-\int_0^v \lambda(\zeta(w))\, dw} \, / \, \zeta(0)=z \biggr)}.$$
En particulier :
$$ \E(g(t,X_t)) = \E(\psi(T_{N_t},Z_t, A_t)).$$
\end{prop}

\begin{dem}
Posons $\bar F_z(v)=\P(T_1>v\, / \, \zeta(0)=z) =\E\left( e^{-\int_0^v \lambda(\zeta(w))\, dw} \, / \, \zeta(0)=z\right).$

Soit $\psi_0$ une fonction mesurable positive définie sur $\R_+ \times \R^d \times \R_+$. On a :
\begin{eqnarray}\label{eqMulti1psi}
\lefteqn{\E(\psi_0(T_{N_t}, Z_t, A_t)\, g(t,X_t)) = }\\
 &  & \sum_{n \ge 0} \E\biggl(1_{\{T_n \le t\}}\, \psi_0(T_n, Y_n, t-T_n) \, \E(g(t,\zeta^{(n+1)}(t-T_n)) \, 1_{\{T_{n+1}-T_n > t-T_n\}} \, / \, Y_n, T_n)\biggr)\nonumber
\end{eqnarray}
Or :
\begin{eqnarray*}
\lefteqn{ \E(g(t,\zeta^{(n+1)}(t-T_n)) \, 1_{\{T_{n+1}-T_n > t-T_n\}} \, / \, Y_n, T_n)}\nonumber\\
& = &  \E\biggl(g(t-T_n+T_n,\zeta^{(n+1)}(t-T_n)) \, M(\zeta^{(n+1)}, \R^d \times ]t-T_n,+\infty[ )  \, / \, Y_n, T_n \biggr)\nonumber\\
& = & \psi(T_n,Y_n, t-T_n)\, \bar F_{Y_n}(t-T_n) = \psi(T_n,Y_n, t-T_n)\, \P(T_{n+1}-T_n > t-T_n\, / \, Y_n, T_n)
\end{eqnarray*}
En reportant dans (\ref{eqMulti1psi}), nous obtenons 
$ \E(\psi_0(T_{N_1},Z_t, A_t)\, g(t,X_t)) = \E(\psi_0(T_{N_t}, Z_t, A_t)$ $ \,\psi(T_{N_t}, Z_t, A_t)).$
\end{dem}

\medskip

On peut  appliquer au SP $(X_t)_{t \ge 0}$ les résultats sur la convergence des processus semi-régénératifs. Notamment, dans le cas non-arithmétique et sous des conditions précisées par exemple dans \cite{Alsmeyer1} ou \cite{Alsmeyer2}, on obtient :
\begin{eqnarray*}
 \E(g(X_t)) &\tend{1cm}_{t \to \infty}& \frac{\displaystyle{\int_{\R^d} \E\left( \int_0^{T_1} g(\zeta(v))\, dv \, / \, \zeta(0)=z)\, dv\right)\, m(dz)}}{\displaystyle{\int_{\R^d} \E(T_1\, / \, \zeta(0)=z)\,m(dz)}}\\
 & = & \frac{\displaystyle{\int_{\R^d} \E \biggl(\int_{\R_+} g(\zeta(v))\, e^{-\int_0^v \lambda(\zeta(w))\,dw}\, dv\, / \, \zeta(0)=z\biggr)\, m(dz)}}{\displaystyle{\int_{\R^d}  \E\biggl(\int_{\R_+} e^{-\int_0^v \lambda(\zeta(w))\,dw}\, dv\, / \, \zeta(0)=z\biggr)\, m(dz)}}
 \end{eqnarray*}
où $m$ est la loi stationnaire de la chaine de Markov $(Y_n)_{n \ge 1}$.

\section{Cas d'un processus intrinsèque markovien}

\begin{theo}
Un processus à changement de rythmes associé à un processus de Markov $\zeta$ et à $M$ donné par (\ref{defM}) est un processus de Markov.
\end{theo}

Ce théorème est une conséquence immédiate du théorème suivant démontré dans \cite{CC}.

\begin{theo}
On suppose que  : 
\begin{itemize}
\item[i.] le processus $\zeta$ est un processus de Markov.
\item[ii.]  pour tout $s \in \R_+$ et toute fonction mesurable positive $\varphi$ définie sur $\R^d \times \R_+$ :
\begin{eqnarray}\label{hyp1}
\lefteqn{\int_{\R^d \times \R_+} \varphi(z,v)\, 1_{\{v>s\}}\, M(\zeta;dz,dv)}\nonumber\\
& \hspace{-1.5cm} = & \hspace{-0.5cm} M(\zeta;\R^d \times ]s,+\infty[)\, \int_{\R^d \times \R_+} \varphi(z,v+s)\,  M(\zeta(s+ \ \cdot \ );dz,dv),
\end{eqnarray}
\item[iii.] pour tout $s \in \R_+$, $M(\zeta;\R^d \times [0,s])$ est mesurable pour la tribu engendrée par  les variables aléatoires $\zeta(v), v \le s$.
\end{itemize}
Alors le SP $(X_t)_{t \ge 0}$ associé à $\zeta$ et $M$ est un processus de Markov.
\end{theo}

\section{Cas d'un processus intrinsèque semi-martingale}

\begin{theo}\label{theoFondSM}
Nous supposons que 
$$ M(\zeta;dz,dv) = \lambda(\zeta(v))\, e^{-\int_0^v \lambda(\zeta(w))\, dx}\, Q(\zeta(v);dz)\, dv$$
 et que $\lambda$ est borné.

Soit $({\cal D}({\cal A}_0), {\cal A}_0)$ un opérateur sur les fonctions à valeurs réelles définies sur $\R_+ \times F$, bornées sur $[0,t] \times F$ pour tout $t>0$. Nous supposons que pour toute fonction $g$ appartenant à ${\cal D}({\cal A}_0)$ :
\begin{enumerate}
\item  le processus $t \to g(t,\zeta(t))$ est une semi-martingale de la forme
 $$g(t,\zeta(t))=g(0,\zeta(0))+ \int_0^t {\cal A}_0g(v, \zeta(v))\, dv + M_t^g$$
où $M_t^g$ est une martingale,
\item la fonction ${\cal A }_0 g$ est bornée sur $[0,t] \times F$ pour tout $t>0$.
\end{enumerate}

Soit $\Psi$ un SP associé à $\zeta$ et $M$ et $g \in {\cal D}({\cal A}_0)$.
Posons :
$$ \tilde {\cal A}_0g(v,z) = {\cal A}_0g(v,z) + \lambda(z) \int_F (g(v,z_1) - g(v,z))\, Q(z;dz_1),$$
On suppose que pour tout $s>0$, $\tau_s g \in {\cal D}({\cal A}_0)$ et ${\cal A}_0 \tau_s g = \tau_s {\cal A}_0 g$. Alors : 
$$ \E(g(t,\Psi_t) ) = \E(g(0, \Psi_0)) + \int_0^t \E(\tilde {\cal A}_0g(s, \Psi_s))\, ds.$$
\end{theo}

\begin{dem}
Nous allons nous appuyer sur la proposition \ref{loiMulti} dont nous reprenons les notations. Nous posons $\E_z( \ \cdot \ ) = \E(\ \cdot \ \, / \, \zeta(0)=z)$.

Remarquons que $\zeta$ étant càd-làg, $\{v : \zeta(v_-) \not= \zeta(v)\}$ est dénombrable donc pour toute fonction mesurable $f : \R_+ \times \R^d \to \R$, nous avons $f(v,\zeta(v_-))\, dv=f(v,\zeta(v))\, dv.$ 

Soit $g$ une fonction bornée appartenant à ${\cal D}({\cal A}_0)$. La formule d'Itô donne 
\begin{eqnarray*}
 e^{-\int_0^t \lambda(\zeta(w))\, dw}\, g(t,\zeta(t)) & = & g(0, \zeta(0)) + \int_0^t e^{-\int_0^s \lambda(\zeta(w)\, dw} {\cal A}_1g(s, \zeta(s))\, ds \\
& & \ + \ \int_0^t e^{-\int_0^s \lambda(\zeta(w)\, dw}\, dM_s^g 
\end{eqnarray*}
avec 
$ {\cal A}_1g(s, x) = {\cal A}_0g(s,x) - \lambda(x)\, g(s,x),$
d'où
\begin{equation}\label{eqMartGen}
 \E_z\left(e^{-\int_0^t \lambda(\zeta(w))\, dw}\, g(t,\zeta(t))\right) - g(0,z) = \E_z\left( \int_0^t  e^{-\int_0^v \lambda(\zeta(w))\, dw} {\cal A}_1 g(v, \zeta(v)) \, dv\right).
\end{equation}

Soit $N$ le noyau du processus de renouvellement markovien $(Y_n, T_n)_{n \ge 0}$.
Nous écrivons $N(x,dy,dv) =1_{\R_+}(v)\, dF_x(v)\, \beta(x,v;dy)$ où $dF_x$ est la loi de $T_1$ sachant $X_0=x$ et $\beta(x,v;dy)$ la loi de $Y_1$ sachant $T_1=v$ et $X_0=x$. D'après (\ref{renouvMulti}) nous avons :
$$1_{\R_+}(v)\, dF_x(v) = \int_{\R^d} N(x,dy,dv) = \E_x( \lambda(\zeta(v))\, e^{-\int_0^v \lambda(\zeta(w))\, dw})\, dv,$$
donc 
$dF_x(v)$ peut s'écrire $dF_x(v) =  \ell(x,v)\, e^{-\int_0^v \ell(x,w)\, dw}\, dv$ et
\begin{equation}\label{ellborne}
\ell(z,v) = \frac{\E_z\left(e^{-\int_0^v \lambda(\zeta(w))\, dw}\, \lambda(\zeta(v))\right)}{\E_z\left(e^{-\int_0^v \lambda(\zeta(w))\, dw}\right)} \le ||\lambda||_\infty.
\end{equation}
Il s'ensuit que $\E(N_t)<+\infty$ (voir par exemple \cite{CC} corollaire 5.9).

Remarquons que 
$\bar F_z(v) = \E_z(e^{-\int_0^v \lambda(\zeta(w))\, dw}) \ge e^{-|| \lambda||_\infty\, v}.$ 
Posons
$$ \psi(s,z,v)=\E_z\left(e^{-\int_0^v\lambda(\zeta(w)\, dw} g(s+v, \zeta(v)) \right)/ \bar F_z(v).$$
La fonction $\psi$ est bornée sur $\R_+ \times \R^d \times [0,t]$ et compte-tenu de  (\ref{eqMartGen}) appliqué à $\tau_s g$, la fonction 
$v \to  \psi(s,z,v)=\E_z\left(e^{-\int_0^v\lambda(\zeta(w)\, dw} g(s+v, \zeta(v))\right)/ \bar F_z(v)$
est absolument continue. Notons $\partial_3 \psi$ "sa" densité, elle est bornée sur $[0,A] \times \R^d \times [0,B]$ pour tous $A>0, B>0$.

Posons :
$$\tilde L \psi(s,z,v) = \partial_3 \psi(s,z,v) + \int_{\R^d} (\psi(s+v,z_1,0) - \psi(s,z,v)) \, \ell(z,v)\, \beta(z,v;dz_1).$$
Le corollaire 5.21 de \cite{CC} donne :
\begin{equation}\label{CKCSMP}
 \E(\psi(T_{N_t}, Z_t, A_t)) =\E(\psi(0,Z_0,0))+ \E\left(\int_0^t \tilde L \psi(T_{N_s}, Z_s, A_s) \right)\, ds.
\end{equation}
La proposition \ref{loiMulti} entraine $ \E(\psi(T_{N_t}, Z_t, A_t))  = \E(g(t,X_t))$. 

Intéressons nous maintenant au deuxième membre de (\ref{CKCSMP}).
Pour  $(s,z) \in \R_+ \times \R^d$, la fonction $v \to \bar F_z(v)\,\tilde L \psi(s,z,v)$ est intégrable sur $[0,t]$.
Une intégration par parties, la relation $dF_z(v) = \bar F_z(v)\, \ell(z,v)\, dv$, l'hypothèse ${\cal A}_0 \tau_s g = \tau_s {\cal A}_0g$ et la formule (\ref{eqMartGen}) appliquée à $\tau_sg$ entrainent :
$$  \int_0^t \bar F_z(v) \tilde L \psi(s,z,v)\, dv = \int_0^t \E_z\left(e^{-\int_0^v \lambda(\zeta(w))\, dw}\, \tilde {\cal A}_0g(s+v, \zeta(v))\right)\, dv.$$
Posons 
$\psi_1(s,z,v) = \E_z\left(e^{-\int_0^v \lambda(\zeta(w))\, dw}\, \tilde {\cal A}_0g(s+v, \zeta(v))\right)/\bar F_z(v).$
Pour tout $t \ge 0$, $\int_0^t \bar F_z(v) \tilde L \psi(s,z,v)\, dv = \int_0^t \bar F_z(v)\, \psi_1(s,z,v)\, dv$ et par conséquent, pour tout $(s,z)\in \R_+ \times \R^d$, $\tilde L \psi(s,z,v)\, dv = \psi_1(s,z,v)\, dv$.
Donc pour $T_n \le t <T_{n+1}$ :
\begin{eqnarray*}
\int_0^t \tilde L \psi(T_{N_s}, Z_s,A_s)\, ds & = &  \sum_{k=0}^{n-1} \int_0^{T_{k+1}-T_k}  \tilde L \psi(T_k, Y_k, v)\, dv + \int_0^{t-T_n} \tilde L \psi(T_n, Y_n, v)\, dv\\
& = & \sum_{k=0}^{n-1} \int_0^{T_{k+1}-T_k}  \psi_1(T_k, Y_k, v)\, dv + \int_0^{t-T_n} \psi_1(T_n, Y_n, v)\, dv\\
& = & \int_0^t \psi_1(T_{N_s}, Z_s,A_s)\, ds.
\end{eqnarray*}
En utilisant à nouveau la proposition \ref{loiMulti} nous obtenons
$ \E\left(\int_0^t \tilde L \psi(T_{N_s}, Z_s, A_s) \right)\, ds =\int_0^t  \E(\tilde {\cal A}_0g(s, X_{s-}))\, ds$, d'où le résultat.
\end{dem}

\begin{cor}\label{corFondSM}
Nous supposons que 
$$ M(\zeta;dz,dv) = \lambda(\zeta(v))\, e^{-\int_0^v \lambda(\zeta(w))\, dx}\, Q(\zeta(v);dz)\, dv$$
 et que $\lambda$ est borné.

Soit $({\cal D}({\cal A}), {\cal A})$ un opérateur sur les fonctions bornées à valeurs réelles définies sur $\R^d$. Nous supposons que pour toute fonction $f$ appartenant à ${\cal D}({\cal A})$ :
\begin{enumerate}
\item  le processus $t \to f(\zeta(t))$ est une semi-martingale de la forme
 $$f(\zeta(t))=f(\zeta(0))+ \int_0^t {\cal A}f(\zeta(v))\, dv + M_t^f$$
où $M_t^f$ est une martingale,
\item la fonction ${\cal A }f$ est bornée sur $\R^d$.
\end{enumerate}

Soit $\Psi$ un SP associé à $\zeta$ et $M$ et $f \in {\cal D}({\cal A})$.
Posons :
$$ \tilde {\cal A}f(z) = {\cal A}f(z) + \lambda(z) \int_{\R^d} (f(z_1) - f(z))\, Q(z;dz_1).$$
 Alors : 
$$ \E(f(\Psi_t) ) = \E(f(\Psi_0)) + \int_0^t \E(\tilde {\cal A} f(\Psi_s))\, ds.$$
\end{cor}

Le résultat suivant est une application du théorème \ref{theoFondSM} lorsque le processus intrinsèque $\zeta$ est un processus de Lévy.

\begin{cor}\label{corLevy} 
Soit $\zeta$ un processus de Lévy d-dimensionnel de mesure de sauts $J$ et de triplet $(\mu, {\cal C},\nu)$, c'est-à-dire
$$ \zeta(t) = \mu\, t + \sqrt{\cal C}\, W_t + \int_0^t \int_{\{x : ||x|| \ge 1\}} x\, J(ds,dx) + \int_0^t \int_{\{x : ||x|| < 1\}} x\,\tilde  J(ds,dx)$$
où $W$ est un brownien standard de dimension $d$ indépendant de $J$, $\cal C$ une matrice de corrélation et $\tilde J(ds,dx) = J(ds,dx) - \nu(dx)\, ds$.

Soit ${\cal D}({\cal A}_0)$ l'ensemble des fonctions $g=g(t,x) : \R_+ \times \R^d \to \R$  bornées de classe $C^{1,2}$, c'est-à-dire continûment différentiable par rapport à la variable temporelle  $t$ et 2 fois continûment différentiable par rapport à la variable spatiale $x$, et dont les dérivées d'ordre 1 et 2 sont bornées sur $[0,t] \times \R^d$. Pour $g \in {\cal D}({\cal A}_0)$, posons
\begin{eqnarray*}
 {\cal A}_0 g(s,x) & = & \frac{\partial g}{\partial s}(s,x) + \sum_{i=1}^d \mu_i \frac{\partial g}{\partial x_i} (s,x)  + \frac{1}{2} \sum_{i=1}^d \sum_{j=1}^d {\cal C}_{i,j}\, \frac{\partial^2 g}{\partial x_i \partial x_j} (s,x) \\ 
   & & \ + \ \int_{\R^d} \biggr(g(s, x + y) - g(s, x) -  \sum_{i=1}^d \frac{\partial g}{\partial x_i} (s, x) \, y_i 1_{\{||y|| < 1\}} \biggr) \, \nu(dy).
   \end{eqnarray*}
   $$ \tilde {\cal A}_0g(s,x) = {\cal A}_0 g(s,x) + \lambda(x) \int_{\R^d} (g(s,y)-g(s,x))\, Q(x;dy).$$

Soit $X$ un SP associé à $\zeta$ et $M$ donné par (\ref{defM}). On suppose que $\lambda$ est borné.

Alors pour toute fonction $g$ appartenant à ${\cal D}({\cal A}_0)$ :
$$ \E(g(t,X_t) ) = \E(g(0, X_0)) + \int_0^t \E(\tilde {\cal A}_0g(s, X_s))\, ds.$$

\end{cor}

\begin{dem}
La formule d'Itô donne :
\begin{eqnarray}\label{LevyA1}
 \lefteqn{ \hspace{-2cm} g(t,\zeta(t)) =   g(0,0) + \int_0^t {\cal A}_0 g(s, \zeta(s))\, ds  + \sum_{i=1}^d  \int_0^t   \frac{\partial g}{\partial x_i} (s, \zeta(s_-))\, \sum_{j=1}^d \sigma_{i,j}\, dW_s^j }\nonumber\\
& & + \int_0^t \int_{\R^d}   \biggr(g(s, \zeta(s_-) + y) - g(s, \zeta(s_-))  \biggr) \,  \tilde J(ds,dy).
  \end{eqnarray}
  D'une part les $\partial g/\partial x_i$ étant bornés, les $\int_0^t (\partial g/\partial x_i)(s, \zeta(s_-))\, dW_s^j$ sont des martingales.
  D'autre part la formule de Taylor entraine
\begin{eqnarray*}
 \lefteqn{\hspace{-3cm} \E \left( \int_0^t | g(s, \zeta(s_-)+y)-g(s, \zeta(s_-) |^2 \nu(dy)\, ds \right)  \le  4 ||g||^2\, t\,  \nu(\{y : ||y||\ge 1\})}\\
& & \  + \ C_t\,  t\,  \int_{\{y : ||y|| <1\}} ||y||^2\, \nu(dy) < +\infty
\end{eqnarray*}
où $C_t$ est une constante qui dépend de $\sup_{1 \le i \le d} \sup_{s \le t, x \in \R^d} | \frac{\partial g}{\partial x_i} (s, x)| $. Par conséquent
$\int_0^t \int_{\R^d}   \biggr(g(s, \zeta_{s-} + y) - g(s, \zeta_{s-})  \biggr) \,  \tilde J(ds,dy)$ est une martingale.

On montre de même que  ${\cal A}_0 g$ est bornée sur $[0,t] \times \R_+$.
\end{dem}

\medskip
En posant $g(s,x) = u(t-s,x)$ on obtient la formule de Feynman-Kac.

\begin{cor}[formule de Feynman-Kac]
Soit $\zeta$ un processus de Lévy d-dimen\-sion\-nel de triplet $(\mu, {\cal C},\nu)$. Soit $X$ un SP associé à $\zeta$ et $M$ donné par (\ref{defM}). On suppose que $\lambda$ est borné.
Soit $u=u(t,x) : \R_+ \times \R^d \to \R$  une fonction  bornée de classe $C^{1,2}$ dont les dérivées d'ordre 1 et 2 sont bornées sur $[0,t] \times \R^d$ pour tout $t$.  On suppose que $u$ vérifie 
$$ \frac{\partial u}{\partial t}(t,x) = {\cal A} u(t,x) + \lambda(x) \int_{\R^d} (u(t,y)-u(t,x))\, Q(x;dy), \quad u(0,x)=h(x),$$
où
\begin{eqnarray*}
 {\cal A} u(t,x)  & = & \sum_{i=1}^d \mu_i \frac{\partial u}{\partial x_i} (t,x)  + \frac{1}{2} \sum_{i=1}^d \sum_{j=1}^d {\cal C}_{i,j}\, \frac{\partial^2 u}{\partial x_i \partial x_j} (t,x) \\ 
   & & \ + \ \int_{\R^d} \biggr(u(t, x + y) - u(t, x) -  \sum_{i=1}^d \frac{\partial u}{\partial x_i} (t, x) \, y_i 1_{\{||y|| < 1\}} \biggr) \, \nu(dy). \end{eqnarray*}
alors $u(t,x) = \E_x(h(X_t))$ pour tout $t$.
\end{cor}

\section{Cas des processus non homogènes en temps}

 Lorsque $\zeta$ est solution d'une équation différentielle stochastique non homogène ou plus généralement un processus d'Itô-Lévy non homogène en temps, on sent bien que le SP $X$ n'est pas défini correctement après un saut $T_n$ $(n \ge 1)$ car il repart comme si celui-ci était l'instant initial. D'ailleurs on ne peut appliquer le théorème  \ref{theoFondSM} car la condition ${\cal A}_0 \tau_sg= \tau_s{\cal A}_0 g $ n'est pas satisfaite. C'est pourquoi nous allons reprendre la définition du SP $X$ dans ce cas et, quitte à être inhomogène, nous allons autoriser la fonction $\lambda$ et le noyau $Q$ à dépendre du temps.

 Soit $\lambda$  une fonction positive définie sur $\R_+ \times \R^d$ et $Q$ est une probabilité de transition de $\R_+ \times \R^d \to {\cal P}(\R^d)$. Posons :
 $$ M(\zeta; dz,dv) = \lambda(v,\zeta(v))\, e^{-\int_0^v \lambda(w, \zeta(w))\, dw}\, Q(v,\zeta(v);dz)\, dv.$$
 
 On construit le SP inhomogène $X$ de la manière suivante. Notons $\zeta^{s,x}=(\zeta^{s,x}(t))_{t \ge 0}$ un processus dont la loi est celle  de $\zeta(s+\ \cdot \  )$ sachant $\zeta(s)=x$. Supposons avoir construit $\zeta^{(1)}, \ldots, \zeta^{(n)}, Y_1, T_1, \ldots, Y_n, T_n$. Soit $\zeta^{(n+1)}$ un processus dont la loi sachant $\zeta^{(1)},\,  \ldots \,  , \zeta^{(n)}, Y_1, T_1, \, \ldots , Y_n=z,   T_n=s$ est la loi de $\zeta^{s,z}$.  La loi de $(Y_{n+1}, T_{n+1}-T_n)$ sachant $\zeta^{(1)}, \ldots, \zeta^{(n+1)}, Y_1, T_1, \ldots, Y_n=z, T_n=s$ est 
$M(\zeta^{(n+1)};dz,dv)$. Pour $T_n \le t < T_{n+1}$ on pose 
$X_t=\zeta^{(n+1)}(t-T_n)$ et $X_{T_{n+1}} = Y_{n+1}$.

\medskip

Si $\xi=(\xi^{(1)}, \xi^{(2)})$ est un processus à valeurs dans $\R_+ \times \R^d$, posons $\tilde M(\xi;ds_1,dz,dv)=M(\xi^{(2)};dz,dv)\, \delta_{\xi^{(1)}(v)}(ds_1)$.
Soit $\tilde \zeta$ un processus dont la loi sachant $\tilde \zeta(0)=(s,x)$ est celle du processus $t \to (s+t, \zeta^{s,x}(t))$.
Le processus $\tilde X=(t,X_t)_{t \ge 0}$ est un SP associé à $\tilde \zeta$ et $\tilde M$ vérifiant $\tilde X^{(1)}=0$.

\medskip
Dans le cas d'un SP inhomogène, $(Y_n, T_n)_{n \ge 1}$ n'est pas un processus de renouvellement markovien. Par contre le processus  $(\tilde Y_n, T_n)_{n \ge 1}$, avec $\tilde Y_n=(T_n, Y_n)$, est un processus de renouvellement markovien de noyau
$$ \tilde N(s,x;ds_1,dx_1,dv)= \E(M(\zeta^{s,x};dx_1,dv))\, \delta_{s+v}(ds_1).$$ 

En appliquant le corollaire \ref{corFondSM} au SP associé à $\tilde \zeta$ et $\tilde M$ on obtient la proposition suivante.

\begin{prop}
On se place dans le cadre de la construction et des notations de ce paragraphe. On suppose que la fonction $\lambda$ est bornée.  

Soit $({\cal D}({\cal A}_0), {\cal A}_0)$ un opérateur sur les fonctions bornées à valeurs réelles définies sur $\R^d$. Nous supposons que pour toute fonction $g$ appartenant à ${\cal D}({\cal A}_0)$ :
\begin{enumerate}
\item  le processus $t \to g(t,\zeta(t))$ est une semi-martingale de la forme
 $$g(t,\zeta(t))=g(0,\zeta(0))+ \int_0^t {\cal A}_0 g(v,\zeta(v))\, dv + M_t^g$$
où $M_t^g$ est une martingale,
\item la fonction ${\cal A}_0 g$ est bornée sur $\R^d$.
\end{enumerate}
 Alors : 
$$ \E(g(t,X_t)) = \E(g(0,X_0)) + \int_0^t \E(\tilde {\cal A}_0 g(s,X_s))\, ds$$
avec :
$$ \tilde {\cal A}_0g(s,x) = {\cal A}_0g(s,x) + \lambda(s,x) \int_{\R^d} (g(s,z) - g(s,x))\, Q(s,x;dz).$$
\end{prop}

\medskip
\noindent
{\bf Exemple.} Lorsque les processus intrinsèques sont des  processus d'Itô-Lévy qui s'é\-cri\-vent 
$$ d\zeta(t) = b(t, \zeta(t_-))\, dt + \sigma(t, \zeta(t_-))\, dW_t + \int a(t, \zeta(t_-),y)\, \tilde J(dt,dy)$$
lorsque $d=1$ (on laisse le lecteur généraliser au cas $d$ quelconque), on a (voir \cite{P}) :
\begin{eqnarray*}
g(t,\zeta(t)) & = & g(0,\zeta(0)) + \int_0^t {\cal A}_0 g(s, \zeta(s)) \, ds + \int_0^t \sigma(s, \zeta(s_-) ) \frac{\partial g}{\partial x} (s, \zeta(s_-))\, dW_s \\
& & + \int_0^t \int (g(s,\zeta(s_-))+ a(s,\zeta(s_-),y)) - f(s, \zeta(s_-)))\, \tilde J(ds,dy)
\end{eqnarray*}
où 
\begin{eqnarray*}
{\cal A}_0g(s,x) & =& \frac{\partial g}{\partial s}(s,x) +  b(s,x) \frac{\partial g}{\partial x} (s,x) + \frac{1}{2} \sigma^2(s,x)\, \frac{\partial^2 g}{\partial x^2}  g(s,x)\\
& & + \int_\R\left(g(s,x+a(s,x,y))-g(s,x) - a(s,x,y) \frac{\partial g}{\partial x}(s,x)\right)\, \nu(dy).
\end{eqnarray*}

\end{document}